\documentclass[]{article}

\usepackage{amsmath,amsfonts,amscd}
\usepackage{xcolor}
\usepackage{graphics}
\usepackage{graphicx}
\usepackage{algorithm, algorithmicx, algpseudocode}
\usepackage{subcaption}
\usepackage{todonotes}
\usepackage{setspace}

\usepackage{hyperref}

\usepackage{graphicx}
\usepackage[margin=2cm]{geometry}

\usepackage{ulem}
\usepackage{todonotes}
\newboolean{SHOWCOMMENTS}
\setboolean{SHOWCOMMENTS}{true}
\ifthenelse{\boolean{SHOWCOMMENTS}}{

}
{

}

\ifthenelse{\boolean{SHOWCOMMENTS}}{
\def\rem#1{{\it\small\color{darkblue} {#1}}}
}
{
\def\rem#1{{}}
}

\ifthenelse{\boolean{SHOWCOMMENTS}}{
\def\SOUT#1{{\sout{#1}}}
}
{
\def\SOUT#1{{}}
}

\ifthenelse{\boolean{SHOWCOMMENTS}}{
\def\SHORTEN#1{{\sout{#1}}}
}
{
\def\SHORTEN#1{{}}
}

\ifthenelse{\boolean{SHOWCOMMENTS}}{

}
{

}

\ifthenelse{\boolean{SHOWCOMMENTS}}{
\setlength{\marginparwidth}{2cm}
\newcommand{\todURi}[1]
{\todo[color=cyan!25,inline]{\footnotesize{\bf Uli:} #1}}
\newcommand{\todUR}[1]%
{\todo[color=cyan!25]{\footnotesize{\bf Uli:} #1} 
}
\newcommand{\todoNK}[1]%
{\todo[color=yellow!25,inline]{\footnotesize{\bf Nils K:} #1}}
\newcommand{\todoDT}[1]%
{\todo[color=cyan!50,inline]{\footnotesize{\bf Dominik T:} #1}}
}
{
\newcommand{\todUR}[1]{}
\newcommand{\todURi}[1]{}
}

\ifthenelse{\boolean{SHOWCOMMENTS}}{
\setlength{\marginparwidth}{2cm}
\newcommand{\todoMSi}[1]{\todo[color=yellow!25,inline]{\footnotesize{\bf ML:} #1}}
\newcommand{\todoMS}[1] {\todo[color=yellow!25]{\footnotesize{\bf ML:} #1}}
}
{
\newcommand{\todoMSi}[1]{}
\newcommand{\todoMS}[1]{}
}

\ifthenelse{\boolean{SHOWCOMMENTS}}{
\setlength{\marginparwidth}{2cm}
\newcommand{\todoCKi}[1]{\todo[color=yellow!25,inline]{\footnotesize{\bf CK:} #1}}
\newcommand{\todoCK}[1] {\todo[color=yellow!25]{\footnotesize{\bf CK:} #1}}
}
{
\newcommand{\todoCKi}[1]{}
\newcommand{\todoCK}[1]{}
}

\newcommand{\Div}{\mbox{div}}

\newcommand{\diag}{\mbox{diag}}

\newcommand{\bv}{{\bf b}}

\newcommand{\dv}{{\bf d}}
\newcommand{\ev}{{\bf e}}
\newcommand{\fv}{{\bf f}}
\newcommand{\gvv}{{\bf g}}

\newcommand{\pv}{{\bf p}}
\newcommand{\qv}{{\bf q}}
\newcommand{\rv}{{\bf r}}

\newcommand{\uv}{{\bf u}}
\newcommand{\vvv}{{\bf v}}
\newcommand{\wv}{{\bf w}}
\newcommand{\xv}{{\bf x}}
\newcommand{\yv}{{\bf y}}
\newcommand{\zv}{{\bf z}}

\newcommand{\Qm}{{\bf Q}}

\newcommand{\Dm}{{\bf D}}

\newcommand{\Rm}{{\bf R}}
\newcommand{\Nm}{{\bf N}}
\newcommand{\Mm}{{\bf M}}
\newcommand{\Am}{{\bf A}}
\newcommand{\Bm}{{\bf B}}

\newcommand{\Id}{{\bf I}}
\newcommand{\Vm}{{\bf V}}

\newcommand{\Wm}{{\bf W}}

\usepackage{authblk}

\title{ Parallel solution of saddle point systems with nested iterative solvers based on the Golub-Kahan Bidiagonalization}

\author{Carola Kruse\thanks{Cerfacs, Toulouse, France},  
Masha Sosonkina\thanks{Old Dominion University, Norfolk, USA},
Mario Arioli \thanks{Libera Universita Mediterranea, Casamassima, Italy},
Nicolas Tardieu \thanks{IMSIA, UMR 9219 EDF-CNRS-CEA-ENSTA, Universit\'e Paris Saclay, Paris, France},
Ulrich R\"ude \thanks{Friedrich-Alexander-Universit\"at Erlangen-N\"urnberg, Erlangen, Germany}}

\begin{document}

\maketitle            

\abstract{We present a scalability study 
of Golub-Kahan
   bidiagonalization for the parallel iterative solution 
   of symmetric indefinite linear systems 
   with a $ 2 \times 2$  block structure. 
   The algorithms have been implemented within the parallel numerical library PETSc. 
   Since a nested inner-outer iteration strategy may be necessary, we investigate different choices for the inner solvers, including parallel sparse direct and multigrid accelerated iterative methods.
   We show the strong and weak scalability of the Golub-Kahan bidiagonalization based iterative method when applied to a two-dimensional Poiseuille flow and to two- and three-dimensional Stokes test problems. }

\setcounter{footnote}{0}
\section{Introduction}\label{sec:intro}

As current and future high-performance computing (HPC) platforms scale,
calculations will be able to increase in size and complexity and take advantage 
of the available processing power and memory.  
At this scale,  HPC applications will increasingly rely on the
parallel numerical libraries
and environments, such as those offered by Trilinos~\cite{trilinos2005} or PETSc~\cite{petsc-web-page} for the solution of
large-scale linear systems using either direct or iterative methods.
Such frameworks help 
to abstract the low-level parallel programming details and enable 
application users and developers to focus on their domain problem at hand.
A report from the U.S. Department of Energy~\cite{dongarra2014applied} notes that ``Numerical
libraries will continue to play an important
role at the exascale'' and that they will allow to share the 
methods implemented therein among ``applications with similar
characteristics''.

In this article, we focus on iterative solvers for indefinite saddle point systems of the type
\begin{align}
 \left(
 \begin{array}{cc}
  \Wm & \Am \\
  \Am^T & 0
 \end{array}
 \right)
  \left(
  \begin{array}{c}
   \wv\\
   \pv
  \end{array}
  \right)
  = \left(
   \begin{array}{c}
   \gvv \\
   \rv
  \end{array}
  \right),\label{eqn:saddlepoint}
\end{align}
with a symmetric positive semi-definite matrix $\Wm\in\mathbb{R}^{m\times m}$ and
$\Am\in\mathbb{R}^{m \times n}$. These systems arise in many applications and
their efficient solution is an active research area. A comprehensive review of
application fields, solvers and preconditioners can be found in
\cite{BeGoLi2005}. Arioli \cite{Ar2013} proposed a new iterative algorithm
by generalizing the standard Golub-Kahan bidiagonalization to matrices of type
(\ref{eqn:saddlepoint}).
In a recent project
jointly with the French electric utility 
EDF, we further investigated this generalized 
Golub-Kahan bidiagonalization solver (GKB) \cite{ArioliKruse2018}.
We have shown 
on the industrial test case of the structural analysis of nuclear reactor
containment buildings that the solver converges in only a
few steps.
A major milestone in the project was the deployment of the developed
GKB-based solver for industrial use at EDF. In several application fields, numerical studies are run in the company with the initially in-house and then later open source finite element software {\sc
code\_aster}\footnote{\url{https://www.code-aster.org}}. 
It is interfaced with PETSc,
which motivated our selection to implement the GKB solver into this framework.
A major advantage of this choice is that the algorithm leverages many PETSc features, such as high-degree of
parallelism and efficient parallel
implementation of basic sparse linear algebra
operations~\cite{petsc-efficient}. 

In this work, we focus on extending and improving our previous research on the parallel performance of the GKB solver \cite{KrSo2020}. By linking with the MKL library for executing dense linear algebra operations, we first improve previously obtained computation times, especially those for MUMPS. As a new test problem, we introduce a 3D Stokes example, which is discretized by Q2-P1 finite elements and we show that the solver is scalable for a fixed problem size. We discuss the weak scalability of the nested inner-outer iterative variants of our solver on the three test cases in two- and three-dimensions. Furthermore, we investigate the portability and present the performance of the algorithm on an AMD architecture.

This paper is organized as follows. Some theoretical aspects of
the GKB algorithm are reviewed in Section~\ref{sec:GGKB}. Then we will comment
on its implementation and usage in PETSc in Section~\ref{sec:petsc}. In
Section~\ref{sec:numerical}, we introduce the Poiseuille flow problem and two- and three-dimensional Stokes test problems and determine required
parameters for the set up and in the stopping criterion. Finally we investigate the strong and weak scalability of the nested iterative solver combinations and show that the GKB method is scalable when an efficient inner solver is used.

\section{Generalized Golub-Kahan Bidiagonalization}\label{sec:GGKB}
We start by summarizing the main results of \cite{Ar2013} which are needed in our further discussion.
The generalized Golub-Kahan bidiagonalization algorithm requires a positive definite (1,1)-block and $\gvv=0$ in the right-hand side.
Depending on the application, $\Wm$ may, however, be only positive semi-definite.
A common method to obtain a positive definite (1,1)-block is to apply the
augmented Lagrangian approach\cite{Ar2013,BeGoLi2005,GoGr2003}. Let $\ker(\Wm)\cap \ker(\Am^T) = \{0 \} $ and $\Nm \in \mathbb{R}^{n\times n}$ be symmetric positive definite.
We modify the upper left block to
\begin{align}
 \Mm := \Wm + \gamma \Am \Nm^{-1} \Am^T \label{eqn:Maug}
\end{align}
for some $0\leq \gamma \leq 1$. 
With the additional transformation 
\begin{align}
\begin{array}{lll}
\uv &= &\wv - \Mm^{-1}(\gvv + \gamma \Am\Nm^{-1}\rv)\\
\bv & = &\rv - \Am^T\Mm^{-1}(\gvv + \gamma \Am \Nm^{-1}\rv),
\end{array}
\label{eqn:trafo_semi_def}
\end{align}
(\ref{eqn:saddlepoint}) is then equivalent to
\begin{align}
\left[
\begin{array}{cc}
\Mm & \Am\\
\Am^T & 0
\end{array}
\right]
\left[
\begin{array}{c}
\uv \\
\pv
\end{array}
\right]
=
\left[
\begin{array}{c}
 0 \\
\bv
\end{array}
\right].\label{eqn:augsys_auglag}
\end{align}
The non-singularity of $\Mm$ follows from $\ker(\Wm)\cap \ker(\Am^T) = \{0\} $. 
This kind of regularization 
can also be applied when $\Wm$ is positive definite, with the goal that for a suitably chosen $\Nm$, we may find that (\ref{eqn:augsys_auglag})
becomes easier to solve than the original system.
In the following, we will use the notation $\Mm$ for a symmetric 
positive definite matrix. We will use $\gamma = 1$ whenever an augmented Lagrangian approach is used and $\gamma=0$ to work with the original matrix $\Wm$ (thus avoiding the matrix transformation in (\ref{eqn:augsys_auglag})). We will not discuss any intermediate value of $\gamma$, as this factor can be included in $\Nm$. Furthermore, we use the Hilbert spaces
\begin{align*}
 \mathcal{M} = \{{\bf v} \in \mathbb{R}^m: \|{\bf v} \|_{\Mm}^2 = {\bf v}^T \Mm {\bf v}\}, \hspace{0.3cm} \mathcal{N} = \{\qv \in \mathbb{R}^n: \|\qv \|_{\Nm}^2 = \qv^T \Nm \qv \}.  
\end{align*}

\subsection{Fundamentals of the Golub-Kahan Bidiagonalization Algorithm}\label{sec:fundGKB}

The (standard) Golub-Kahan bidiagonalization procedure has been widely used in the computation of the singular value decomposition of rectangular matrices. Let $\tilde{\Am}\in\mathbb{R}^{m x n}$, then we search for two unitary matrices $\tilde{\Qm}^{n \times n}$ and $\tilde{\Vm}^{m \times m}$, such that $\tilde{\Vm}^T \tilde{\Am} \tilde{\Qm} = \Bm$, where
\begin{align*}
\Bm =
\left[ \begin{array}{ccccc}
\alpha_1 & \beta_2 &  0 & \cdots & 0 \\
0 & \alpha_2 & \beta_3 &   \ddots & 0 \\
\vdots &\ddots  & \ddots  & \ddots &\ddots  \\
0 & \cdots & 0 &\alpha_{n-1} & \beta_{n}   \\
 0 & \cdots &  0 &  0 & \alpha_n
\end{array}\right]  .
\end{align*}
In \cite{GoKa1965,PaSa1982,Sa1995}, several algorithms for the bidiagonalization are presented that can be applied to $\tilde{\Am}$. Here, we  will specifically analyze one of the variants known as the Craig-variant \cite{PaSa1982,Sa1995,Sa1997}. With the transformations $\Am = \Mm^{1/2}\tilde{\Am}\Nm^{1/2}$, $\Qm=\Nm^{-1}\tilde{\Qm}$  and $\Vm=\Mm^{-1}\tilde{\Vm}$, the Golub-Kahan bidiagonalization can be generalized into seeking the matrices $\Qm, \Vm$ and $\Bm$, such that 
\begin{align}
\left\lbrace
\begin{array}{r@{}c@{}ll@{}l}
\Am \Qm &=& \Mm \Vm \left[ \begin{array}{c}\Bm\\ 0\end{array}  \right], &\qquad \Vm^T \Mm \Vm &= \Id_m \\
&&\\
\Am^T \Vm &=& \Nm \Qm \left[  \Bm^T  ; 0 \right], &\qquad \Qm^T \Nm \Qm &= \Id_n
\end{array}
.\right. \label{eqn:GKalg}
\end{align}
For a more detailed derivation, we refer to \cite{Ar2013,OrAr2017}. By the change of variables $\uv := \Vm \hat{\zv}$ and $\pv := \Qm \hat{\yv}$ and by multiplying the system from the left by the block diagonal matrix blockdiag($\Vm^T,\Qm^T$), the augmented system can be transformed with (\ref{eqn:GKalg}) into
\begin{eqnarray*}
\left[ \begin{array}{ccc}
\Id_n   & 0              & \Bm\\
0         & \Id_{m-n} & 0 \\
\Bm^T & 0              & 0
\end{array}\right]
\left[ \begin{array}{c}
\hat{\zv}_1 \\ \hat{\zv}_2 \\ \hat{\yv}
\end{array}\right]  =
\left[ \begin{array}{c}
0 \\ 0 \\ \Qm^T \bv
\end{array} \right] .
\end{eqnarray*}
It follows immediately that $\hat{\zv}^T = (\hat{\zv}^T_1, \hat{\zv}^T_2) = (\hat{\zv}^T_1, 0)$. Consequently, $\uv$ only depends on the first $n$ columns of $\Vm$  and, thus, the system reduces to
\begin{align*}
\left[ \begin{array}{cc}
\Id_n   & \Bm\\     \Bm^T       & 0
\end{array}\right]
\left[ \begin{array}{c}
\hat{\zv}_1  \\ \hat{\yv}
\end{array}\right]  =
\left[ \begin{array}{c}
 0 \\ \Qm^T \bv
\end{array} \right] .
\end{align*}
The GKB algorithm can be set such that $\Qm^T \bv = \|\bv\|_\Nm \ev_1$ by choosing $\qv_1 = \Nm^{-1} \bv / \|\bv\|_{\Nm^{-1}}$.
We denote the iterates of $\hat{\zv}_1$ by $\zv_k$ and the entries of $\zv_k$ by $\zeta_j$, $j=1,..,k$, i.e. $\zv_k^T = (\zeta_1,..,\zeta_k)$. 
In \cite{Ar2013}, it is proved that by taking 
advantage of the recursive properties of the standard Golub-Kahan algorithm \cite{GoKa1965}, and using some of
the results of \cite{PaSa1982}, we can obtain the fully  recursive Craig's variant algorithm (Algorithm~\ref{alg:GKB}).
We highlight that in each iteration two linear systems, one for $\Mm$ and one for $\Nm$ must be solved. In the following, we use exclusively $\Nm=\frac{1}{\nu}\Id$, so that the inversion of $\Nm$ has only negligible cost. On the other hand, when applying the augmented Lagrangian approach (\ref{eqn:augsys_auglag}), the matrix $\Mm$ usually suffers of ill-conditioning for increasing values of $\nu$. The inversion of $\Mm$ might then become costly, whereas the GKB algorithm converges usually in fewer iterations \cite{ArioliKruse2018}.

\begin{algorithm}
  \caption{Golub-Kahan bidiagonalization algorithm}
  \label{alg:GKB}
  \begin{algorithmic}
  \Require{$\Mm , \Am , \Nm, \bv$, maxit}
  \State{$k=0$; $\beta_1 = \|\bv\|_{\Nm^{-1}}$;  $\qv_1 = \Nm^{-1} \bv / \beta_1$}
  \State{$\wv = \Mm^{-1} \Am \qv_1$; $\alpha_1 = \|\wv\|_{\Mm}$; $\vvv_1 = \wv / \alpha_1$}
  \State{$\zeta_1 = \beta_1 / \alpha_1$; $\dv_1=\qv_1/ \alpha_1$; $\pv^{(1)} = - \zeta_1 \dv_1$; $\uv_{1} = \zeta_1 \vvv_1$;}
  \While{convergence = false  and $k < $ maxit}
  \State{$k = k + 1$}
  \State{$\gvv = \Nm^{-1} \left( \Am^T \vvv_k - \alpha_k \Nm \qv_k  \right) $; $\beta_{k+1} = \|\gvv\|_{\Nm}$}
  \State{$\qv_{k+1} = \gvv / {\beta_{k+1}}$}
  \State{$\wv = \Mm^{-1} \left(  \Am \qv_{k+1} - \beta_{k+1} \Mm \vvv_{k} \right)$; $\alpha_{k+1} = \|\wv\|_{\Mm}$}
  \State{$\vvv_{k+1} = \wv / {\alpha_{k+1} }$}
  \State{$\zeta_{k+1} = - \dfrac{\beta_{k+1}}{\alpha_{k+1}} \zeta_k$}
  \State{$\dv_{k+1} = \left( \qv_{k+1} - \beta_{k+1} \dv_k \right) / \alpha_{k+1} $}
  \State{$\uv_{k+1} = \uv_{k} + \zeta_{k+1} \vvv_{k+1}$; $\pv_{k+1} = \pv_{k} - \zeta_{k+1} \dv_{k+1}$}
  \State{$\left[ \right.  $ convergence $\left. \right] $ = check$(\zv_k, \dots)$}
  \EndWhile \\
  \Return $\uv_{k+1}, \pv_{k+1}$
  \end{algorithmic}
\end{algorithm}

\subsection{Stopping Criterion}\label{sec:stop}
The statement 'check($\zv_k$)' in Algorithm~\ref{alg:GKB} is yet undefined. In this section, we review a lower bound estimate of the error in energy norm as stopping criterion that was initially proposed by Arioli \cite{Ar2013}. The error $\ev_{k} = \uv - \uv_{k}$ can be expressed, using the relations in (\ref{eqn:GKalg}) and the \Mm-orthogonality property of $\Vm$, by
\begin{align*}
\| \ev_{k} \|_{\Mm}^2 = \sum_{j=k+1}^n \zeta_j^2 = \Big|\Big| \hat{\zv} - \left[ \begin{array}{c}
\zv_k \\ 0
\end{array}\right] \Big|\Big|_2^2.
\end{align*} 
To compute $\ev_{k}$, we thus need $\zeta_{k+1}$ to $\zeta_{n}$, which are available only after the full $n$ iterations of the algorithm. Given a threshold $\tau < 1$ and an integer $d$, we define lower bounds of $\| \ev_{k} \|_{\Mm}^2$  and  $\|\uv\|_{\Mm}$ by 
\begin{align*} 
\xi_{k,d}^2 &= \sum_{j=k+1}^{k+d+1} \zeta_j^2 < \| \ev_{k} \|_{\Mm}^2 , \quad 
\sum_{j=1}^{k+d+1} \zeta_j^2 < \|\uv \|_\Mm^2 
\end{align*}
and, then, by them a stopping criterion  
\begin{align}
&\mbox{\textbf{if}~~~} \xi_{k,d}^2 \leq \tau \sum_{j=1}^{k+d+1} \zeta_j^2 \mbox{~~~~\textbf{stop}}.
 \label{eqn:stop}
\end{align}
$\xi_{k,d}$ measures the error at step $k-d$, but as the following $\uv_{k}$ minimize the error due to an important property of minimization of Craig's algorithm \cite{Sa1995}, we can safely use the last ones. 
For computational examples underlying the efficiency of this lower bound stopping criterion, we refer the reader to \cite{Ar2013,ArioliKruse2018}. In terms of numerical cost, this lower bound estimate is very inexpensive to compute.

\section{Implementation and Usage Details}\label{sec:petsc}

We have implemented the GKB solver in the PETSc {\tt PCFIELDSPLIT} \cite[Chapter 4.5]{petsc-user-ref} 
environment and it is available in the 3.11 release. {\tt PCFIELDSPLIT} provides preconditioners and solvers for block-matrices, as for example several variants of Schur complement preconditioners for 2x2 block matrices. 

Similar to many solution methods in PETSc, the GKB method can be used as either a preconditioner or solver. To obtain GKB as solver, the standard PETSc options are to be set as 
\textit{-ksp\_type preonly -pc\_type fieldsplit -pc\_fieldsplit\_type gkb}.
The solver can only be used for symmetric block matrix systems with zero (2,2)-block
as in Eq.~(\ref{eqn:saddlepoint}). If the matrix is not symmetric, the code
will stop with an error message. The (1,1)-block may be positive semi-definite
or definite. An augmented Lagrangian
approach must be used to ensure the non-singularity of the matrix in the first
case and it can be used to obtain a potentially better convergence in
the second case (see Section~\ref{sec:GGKB}). The GKB PETSc options are 
\begin{center}
\begin{tabular}{|p{3.9cm}|p{8.1cm}|}
\hline
\textit{-pc\_fieldsplit\_gkb\_nu} & $\nu > 0$: Eq.
(\ref{eqn:Maug}) is used with $\gamma=1, \Nm=\frac{1}{\nu}\Id$. \\
&$\nu=0$:
Original system is used, i.e. $\gamma=0$ and $\Nm=\Id$.\\
\hline
\textit{-pc\_fieldsplit\_gkb\_delay} & The delay $d$ in the lower bound stopping criterion of Section
\ref{sec:stop}\\
\hline
\textit{-pc\_fieldsplit\_gkb\_tol} & Stopping tolerance $\tau$ of the solver. \\
\hline
\textit{-pc\_fieldsplit\_gkb\_maxit} & Maximal number of iterations.\\
\hline
\textit{-pc\_fieldsplit\_gkb\_monitor} & Displays the lower bound estimate at each iteration.\\
\hline
\end{tabular}
\end{center}
In general, $\Nm$ may be any kind of positive definite matrix.
In our PETSc implementation, the matrix is however restricted to
$\Nm = \frac{1}{\nu}\Id$. 
The augmented Lagrangian approach is switched off with \textit{-pc\_fieldsplit\_gkb\_nu 0}.
This is done for the convenience of not passing the
parameter $\gamma$, but corresponds to $\gamma=0$, $\Nm = \Id$.

A considerable advantage of the integration of the GKB iterative method
in PETSc is the availability of a large choice of solver-preconditioner combinations for the inner solution step of linear systems
of type $\Mm \xv = \fv$ in Algorithm~\ref{alg:GKB}.
Although the outer loop in the GKB method is sequential, each matrix or vector operation is fully parallel and scalability is achieved by the inner solvers (see Section~\ref{sec:numerical}). 

\section{Numerical Experiments} \label{sec:numerical}

Iterative solvers for the Stokes system have been a field of intensive research. These include preconditioned Krylov subspace methods \cite{elman2002performance,elman2014finite,ur2011iterative}
and multigrid methods \cite{brandt2011multigrid}. Parallel multigrid methods for the Stokes system have been studied in \cite{gmeiner2016quantitative,gmeiner2015towards}. 
In this section, we will apply the the GKB iterative solver to two discretizations of the Stokes equations in two and three dimensions. A comparative study between the GKB and the previously cited methods would however be out of the scope of this paper. Here, we will focus on comparing its parallel performance for different inner solvers as well as to the sparse direct solver MUMPS \cite{MUMPS:1} applied to the overall system. In particular, we discuss choices for $\nu$ as well as the stopping tolerances $\tau$ of the GKB method and  $\tau_{in}$ of its inner iterative solvers.

\paragraph{Experimental Set-up} The calculations are executed on the cluster Kraken of the Cerfacs  computing resources. 
The Kraken cluster includes 121 compute nodes equipped with two Intel Skylake processors at 2.3 Ghz, each of them has 18 cores and share 96 GB DDR4 memory.
Unless otherwise stated, we use a power of two number of cores. In the numerical computations, we first fill up one node up to $2^5$ cores, and for higher counts  32 out of the 36 cores are used per node.
The MPI tasks are evenly distributed among the processors of a node. 
For the largest case {\tt Prob 4} ($>25\cdot 10^{6}$ unknowns) in Section~\ref{sec:discerror}, the computations with MUMPS are done on 32 of 36 cores on one ``fat'' compute node with 768 GB of memory. 
This is necessary as MUMPS exits with a memory error on the standard compute nodes. PETSc and the examples are compiled with gcc 8.3 and openmpi 4.0.1. For linear algebra kernels, MKL version 2018.1.163 is used. The computation times (in seconds) are obtained with the PETSc profiling option \textit{log\_view}, from which we present the time for \textit{ksp\_solve}.

\subsection{Poieseuille Flow}
We consider a viscous, laminar flow in a 2D channel $\Omega=[0,2]\times [0,1]$ with parabolic velocity profile and linear pressure drop, i.e.
\begin{align*}
\uv(x,y)&=(4y(1-y), 0)\\
\pv(x,y) &= 8(2-x).
\end{align*}
This Poiseuille flow problem is the exact solution of the 2D Stokes problem 
 \begin{align}
 \begin{array}{rl}
  -\Delta \uv + \nabla \pv &= 0, \\
  \Div(\uv) &= 0
  \end{array}
  \label{eqn:stokes}
 \end{align}
with no-slip boundary conditions.
We adapt
{\tt ex70.c}\footnote{\url{https://www.mcs.anl.gov/petsc/petsc-dev/src/snes/examples/tutorials/ex70.c.html}} 
given in the PETSc {\tt SNES} section for the simulations.  
The domain $\Omega$ is discretized into a Cartesian grid, using $n_x, n_y$ elements in $x$- and $y$-direction, respectively. The equations are approximated by a cell-centered
 co-located finite volume method,
 where, after application of Gauss's divergence theorem,
 the gradient term is discretized by central differencing and linear interpolation is used for $\pv$ and $\uv$ in the momentum equation \cite{ferziger2001,Klaij2015}.
 This discretization leads to a block system of the form  (\ref{eqn:augsys_auglag}) with a symmetric positive definite matrix $\Mm$. Let, in Matlab notation, $\Dm = \diag(\diag(\Mm))$ and $\Rm=\diag(\diag(\Am^T \Dm^{-1} \Am))$.
 To equilibrate the different blocks of (\ref{eqn:augsys_auglag}), we scale the system from the left and the right by the block diagonal
matrix $\mbox{blockdiag}(\Dm^{-1/2}, \Rm^{-1/2})$.

\subsubsection{Discretization Error}\label{sec:discerror}
In a first experiment, we determine the discretization errors of the model for a sequence of mesh sizes.
This will indicate the necessary stopping tolerance for the iterative algorithms, which thus ensures a fair comparison between the direct and iterative solvers.
Once the accuracy of the iterative solution 
falls below the interpolation error on the nodes ($O(h^2)$ for FEM), 
the solver can stop, as no further significant improvement of the solution accuracy on the exact solution of the PDE can be achieved. 
To get an accurate estimate of the discretization error, we do not take advantage of the block structure of the matrix and solve the complete system with MUMPS \cite{MUMPS:1}. In Table \ref{tab:mumps1}, 
we present the number of degrees of freedom for the (1,1)-block $\Mm$ ($\mbox{dof}_{\Mm}$), the number of constraints ($\mbox{dof}_{\Am}$) and the discretization errors of $\uv_h$ and $\pv_h$ in the 2- and energy-norms
\begin{align*}
 \mbox{err}_2^{\uv} = \frac{1}{n_x n_y}\|\uv_h-\uv\|_2, \hspace{0.4cm}  \mbox{err}^{\pv}_2=\frac{1}{n_x n_y}\|\pv_h-\pv\|_2, \hspace{0.4cm} \mbox{err}_{\Mm}^u=\frac{\|\uv_h-\uv\|_{\Mm}}{\|\uv\|_{\Mm}}.
\end{align*}

\begin{table}
\centering
\caption{Discretization information for the model problem, solved with MUMPS.}
\begin{tabular}{l|r|r|r|r|r|r|r|r} 
name& \multicolumn{1}{|c|}{$n_x$} & \multicolumn{1}{|c|}{$n_y$} & \multicolumn{1}{|c|}{dof} & \multicolumn{1}{|c|}{$\mbox{dof}_{\Mm}$} & \multicolumn{1}{|c|}{$\mbox{dof}_{\Am}$}& \multicolumn{1}{|c|}{$\mbox{err}_2^{\uv}$} & \multicolumn{1}{|c|}{$\mbox{err}_2^{\pv}$}& \multicolumn{1}{|c}{$\mbox{err}_{\Mm}^{\uv}$}  \\
 \hline
{\tt Prob 1}& 512 & 256  & 393\;216   & 262\;144   & 131\;072  & 6.50e-06 & 1.56e-02 & 4.01e-05   \\
{\tt Prob 2}& 1024& 512  & 1\;572\;864  & 1\;048\;576  & 524\;288  & 1.63e-06 & 7.81e-03 & 1.10e-05  \\
{\tt Prob 3}& 2048& 1024 & 6\;291\;456  & 4\;194\;304  & 2\;097\;152 & 4.06e-07 &  3.90e-03 & 2.57e-06  \\
{\tt Prob 4}& 4096& 2048 & 25\;165\;824 & 16\;777\;216 & 8 388\;608 & 1.02e-07 &  1.95e-03 & 6.45e-07   
 \end{tabular}
\label{tab:mumps1}
\end{table}

\subsubsection{GKB Algorithm -- Direct Inner Solver}
We next discuss the GKB algorithm with the direct inner solver MUMPS.
Although $\Mm$ is symmetric positive definite in our application,
we apply the augmented Lagrangian approach (\ref{eqn:Maug})--(\ref{eqn:augsys_auglag}) with $\Nm = \frac{1}{\nu}\Id$ and show the influence of $\nu$ on the number of GKB iterations and the computation time.
We also present results for the case without augmented Lagrangian approach (i.e. $\gamma=0$ and $\Nm=\Id$).
Since the matrices are all symmetric positive definite,
we use the Cholesky factorization in MUMPS,
which switches off the pivoting and leads to a better performance.
The stopping tolerance of the GKB method $\tau$ is chosen as $\tau = 1/n_y^2 \approx 10^{-6}$  such that we can have \textit{superconvergence}
at the mesh nodes and, thus, a smooth reconstruction of the solutions.
Indeed, $\tau$ is of the same order of magnitude as the energy discretization error $\mbox{err}_{\Mm}^{\uv}$ in Table~\ref{tab:mumps1}.
The delay in the lower bound stopping criterion is chosen as $d=5$.
Results for {\tt Prob 3} and the choice of $\nu$ are presented in Table~\ref{tab:eta_gkb1e-7_20481024}.
Of the $\nu$-values tested, the fastest simulation is obtained for $\nu=10$.
This will thus be our choice in the following experiments.
For a discussion about the influence of the parameter $\nu$ on the convergence of the algorithm on different problem settings, we refer to \cite{Ar2013,ArioliKruse2018}.

\begin{table}
\caption{Choice of $\Nm = \frac{1}{\nu}$ for {\tt Prob 3}. $\tau=10^{-6}$, $d=5$ and 32 cores}
\begin{center}
 \begin{tabular}{r|r|r|r|r|r|r}
 \multicolumn{1}{c|}{$\nu$} & \multicolumn{1}{|c|}{$\mbox{err}_2^{\uv}$} & \multicolumn{1}{|c|}{$\mbox{err}_2^{\pv}$} & \multicolumn{1}{|c|}{$\mbox{err}_{\Mm}^{\uv}$} & \multicolumn{1}{|c|}{l.b. estimate}  & \multicolumn{1}{|c|}{GKB iter} & \multicolumn{1}{|c}{time, $s$}\\
 \hline
0   & 4.20e-07 &  5.42e-04 & 2.59e-06 &  9.87e-07 &106 & 89\\
1   & 4.10e-07 &  5.56e-04 & 2.59e-06 &  9.84e-07 &59  & 75\\
10  & 4.06e-07 &  7.58e-04 & 2.57e-06 &  9.92e-07 &42  & 70 \\
100 & 4.06e-07 &  1.60e-03 & 2.57e-06 &  9.83e-07 &53   & 73 \\
\end{tabular}
\end{center}
\label{tab:eta_gkb1e-7_20481024}
\end{table}

\subsubsection{GKB Algorithm -- Iterative Inner Solver}\label{sec:initer}
We will study inner iterative solvers for the solution of the linear systems involving $\Mm$ in Algorithm~\ref{alg:GKB}. This is necessary when $\Mm$ is too large to be solved with a direct method or advantageous if it contains a structure favorable for highly scalable iterative solvers (e.g., multigrid).
Since in our example $\Mm$ is the stiffness matrix for the Laplacian,
we decided to use CG and flexible gmres (denoted \textit{fgmres}, which allows any iterative solver as a preconditioner) preconditioned by BoomerAMG of the library {\tt Hypre} in PETSc. 
For the preconditioning step, we applied one V-cycle with a symmetric hybrid SOR/Jacobi relaxation scheme.
After numerical tests, also trying different numbers of V- and W-cycles, we report that both variants and BoomerAMG do not converge when applied to the augmented Lagrangian matrix (\ref{eqn:augsys_auglag}).
We, thus, use $\gamma=0$ and $\Nm = \Id$.
We found also that the tolerance $\tau_{in}$ of the inner iterative solver has to be set no less than one order of magnitude 
lower than the outer one to obtain a solution of required accuracy.
It furthermore depends on the employed stopping criteria in the implementation of the inner iterative solver and its compatibility with (\ref{eqn:stop}), which results here into different tolerances for CG and fgmres. 
A more precise analysis is beyond the scope of this paper and is left for future work.
Results are presented in Table~\ref{tab:gkb_innerouter} where the errors given are those of the CG method. The errors for fgmres are equivalent and of sufficient accuracy compared to Table~\ref{tab:mumps1}.
Note that the fgmres method results in a lower computation time than CG, although the algorithm is more complex.
In our investigation of this matter, we have noticed that at each GKB iteration, the inner CG solver needs on average one iteration more than fgmres and most of the time is spent in the setup of the AMG solver. 
\begin{table}
\centering
\caption{Inner-outer GKB algorithm with $d=5$ and CG/fgmres-BoomerAMG on 32 cores.}
\begin{tabular}{l|r|r|r|l|l|l|r|r|r} 
  & \multicolumn{3}{|c|}{tolerances}   &   \multicolumn{3}{|c|}{}& \multicolumn{1}{|c|}{Iter}  & \multicolumn{2}{|c}{time, s: GKB with} \\
  \hline
\multicolumn{1}{c|}{name} &  \multicolumn{1}{c|}{$\tau$} & \multicolumn{1}{c|}{CG} &\multicolumn{1}{c|}{fgmres} & \multicolumn{1}{|c|}{$\mbox{err}_2^{\uv}$} & \multicolumn{1}{|c|}{$\mbox{err}_2^{\pv}$} & \multicolumn{1}{|c|}{$\mbox{err}_{\Mm}^{\uv}$} & \multicolumn{1}{|c|}{GKB} & \multicolumn{1}{|c|}{CG} & \multicolumn{1}{|c}{fgmres}\\
 \hline
{\tt  Prob 1} & 1e-05 & 1e-06 & 1e-07 & 6.53e-06 &  2.57e-03 &  4.04e-05 &  60 & 3  &  3 \\
{\tt Prob 2} & 1e-06 & 1e-07 & 1e-08 & 1.62e-06 &  1.55e-03 &  1.02e-05 &  173 & 41 &  34 \\
{\tt Prob 3} & 1e-06 & 1e-07 & 1e-08 & 4.16e-07 &  5.41e-04 &  2.59e-06 &  106 & 115  & 100 \\
{\tt Prob 4} & 1e-07 & 1e-08 & 1e-09 & 1.03e-07 &  2.75e-04 &  6.47e-07 & 220& 1685 & 1424 \\
 \end{tabular}
\label{tab:gkb_innerouter}
\end{table}

\subsubsection{Prob 3: Strong scaling}\label{sec:parallelProb3}
\begin{figure}
\centering
\begin{subfigure}{0.45\textwidth}
\includegraphics[width=\textwidth]{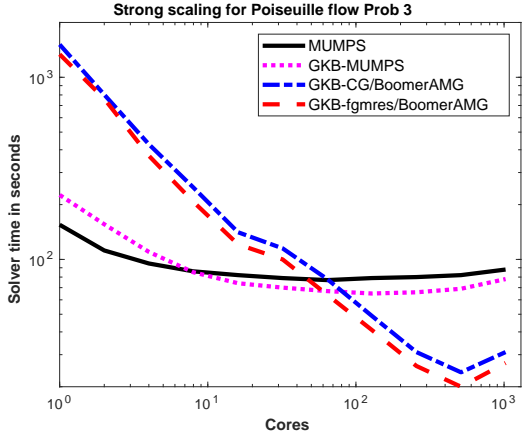}
\caption{Strong scaling Poiseuille flow Prob 3.}
\label{fig:erru_p}
\end{subfigure}
\begin{subfigure}{0.45\textwidth}
\includegraphics[width=\textwidth]{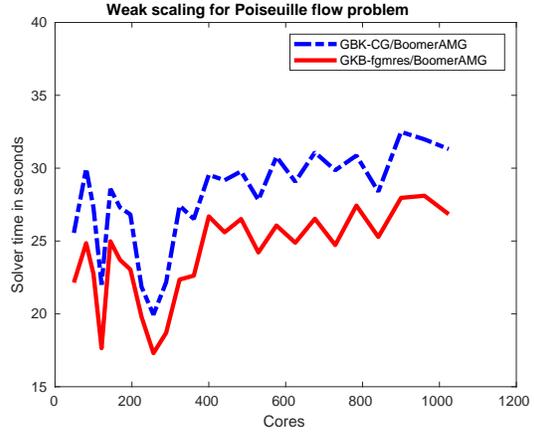}
\caption{Weak scaling Poiseuille flow Prob 3.}
\label{fig:weak_prob3}
\end{subfigure}
\caption{Scaling behavior for the Poiseuille flow problem on Kraken.}
\end{figure}

Strong scaling results for {\tt Prob 3} and the four previously discussed methods are presented in Fig.~\ref{fig:erru_p} and Table~\ref{tab:scalePoiseuille}.
The stopping tolerance for the outer GKB iteration is $\tau=10^{-6}$ and the stopping tolerances for the inner iterative solvers are given in Table~\ref{tab:gkb_innerouter}. 
Furthermore, we choose $\nu=10$ for GKB-MUMPS and the delay $d=5$ in (\ref{eqn:stop}). We observe that on one core, MUMPS applied to the overall system is by almost one order of magnitude faster than GKB-CG and about 9 times faster than GKB-fgmres. 
Once the computations 
use more than 8 cores, all the three GKB variants are faster than MUMPS
in standalone mode.
Among them, the two inner-outer iterative methods with either CG or fgmres clearly outperform GKB-MUMPS from 64 cores onwards, which starts to level off at about 64 cores. 
The performance of the iterative variants starts to saturate at the 1024 core count, thereby indicating the predominance of communication over computation with local, per core, subproblems becoming too small. 
At 512 cores, however, both of them show a speed-up of 63 or more. We emphasize here that CG/fgmres takes advantage of the multigrid preconditioner, which is known to scale well for discretizations of elliptic partial differential equations and, in particular, the Laplace operator. 
Although, as noted earlier, the fgmres method performs better than CG does, the scaling behavior of the two methods is similar. 

 \begin{table}
 \centering
 \caption{Solver time and strong scaling for {\tt Prob 3}. } 
 \begin{tabular}{l||l|l||l|l|p{1.6cm}|p{1.6cm}|p{1.6cm}|p{1.5cm}}
  &\multicolumn{2}{c||}{MUMPS}  & \multicolumn{2}{c|}{GKB-MUMPS} & \multicolumn{2}{c|}{GKB-CG} &\multicolumn{2}{c}{GKB-fgmres}\\
 & \multicolumn{2}{c||}{}&  \multicolumn{2}{c|}{\scriptsize{($\nu$=10)}} & \multicolumn{2}{c|}{\scriptsize{($\nu$=0,$\tau$=$10^{-6}$,$\tau_{in}$=$10^{-7}$)}} &\multicolumn{2}{c}{\scriptsize{($\nu$=$0$,$\tau$=$10^{-6}$,$\tau_{in}$=$10^{-8}$)}}\\
 cores & time (s) & scale & time (s) & scale & time (s) & scale& time (s)& scale\\ 
  \hline
 1 & 155 & & 226 & & 1513 & & 1337 &  \\
 2 & 112 & 1.4 &  156 & 1.5 & 803 & 1.9 & 762 & 1.8   \\
 4 & 95 & 1.6 &  110 & 2.0 & 429 &3.5 & 371 & 3.6  \\
 8 & 86 & 1.8 &  85 & 2.7 & 248 & 6.1 & 207 &  6.5 \\
 16 & 82 & 1.9 & 74 & 3.1& 141 &10.7 & 121 &  11.0 \\
 32 & 79 & 2.0 & 70 & 3.2& 115 & 13.2& 100 &  13.4 \\
 64 & 77 & 2.0 & 67 & 3.4& 78 & 19.4 & 64 &  20.9 \\
 128 & 79  & 2.0 &65 & 3.5& 49 & 30.9 & 41 &  32.6 \\
 256 & 80 & 1.9 &66 & 3.4& 31 & 48.8 & 26 &  51.4 \\
 512 &82 & 1.9 &69 & 3.3& 24 & 63.0 & 20 & 68.9 \\
 1024 &88 & 1.8& 78& 2.9 &31 & 48.8 & 27 & 49.5
  \end{tabular}
  \label{tab:scalePoiseuille}
 \end{table}

\subsubsection{Prob 4: Strong scaling}\label{sec:parallelProb4}
\begin{figure}
    \centering
      \includegraphics[width=8cm]{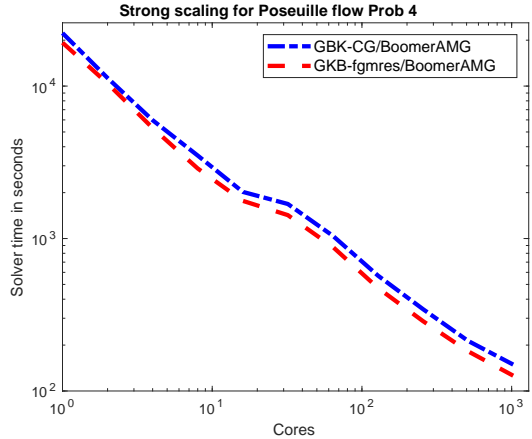}
      \captionof{figure}{Strong scaling Poiseuille flow Prob 4.}
      \label{fig:str_70_4096}
\end{figure}

We repeat the tests of the previous section for the largest case, {\tt Prob 4}, having $16.7\cdot 10^6$ unknowns in the (1,1)-block and $8.3\cdot 10^6$ constraints. On our local machine Kraken, MUMPS requires more memory than available, either for the solution of the entire system ($25\cdot 10^6$ unknowns) or for the GKB-MUMPS solver. We will thus only present the strong scaling performance of GKB-CG/BoomerAMG and GKB-fgmres/BoomerAMG in Figure~\ref{fig:str_70_4096} and Table~\ref{tab:scalePoiseuille4096}. As indicated in Table \ref{tab:gkb_innerouter}, we use as stopping tolerances $\tau =10^{-7}$ for the GKB outer iteration and $\tau_{in} = \{10^{-8},10^{-9}\}$ for the CG and fgmres inner solver, respectively. Furthermore, we use no augmented Lagrangian approach ($\nu=0$) and $d=5$ in the stopping criterion (\ref{eqn:stop}). As before, we obtain faster computations in absolute time with fgmres. Furthermore, due to a much larger size of {\tt Prob 4}, the performances in terms of scalability of both methods are better than those for {\tt Prob 3} are so. GKB-CG shows a speedup of 105 at 512 cores (compared to 63 above) and of about 1.4 when increasing the number of cores from 512 to 1024, reaching a total of 150. We observe the same behavior for GKB-fgmres with 106 at 512 cores (instead of $\sim$69) and a further speedup of 1.4 from 512 to 1024 cores, for a total speedup of 153.

\begin{table}
 \caption{Solver time and strong scaling for {\tt Prob 4}\label{tab:scalePoiseuille4096}}
    \centering
    \begin{tabular}{l||p{1.3cm}|p{1.3cm}|p{1.3cm}|p{1.3cm}}
  & \multicolumn{2}{c|}{GKB-CG} &\multicolumn{2}{c}{GKB-fgmres}\\
 & \multicolumn{2}{c|}{\scriptsize{($\nu$=0,$\tau$=$10^{-7}$,$\tau_{in}$=$10^{-8}$)}} &\multicolumn{2}{c}{\scriptsize{($\nu$=$0$,$\tau$=$10^{-7}$,$\tau_{in}$=$10^{-9}$)}}\\
 cores & time (s) & scale& time (s)& scale\\ 
  \hline
 1 &22271 & & 19248& \\ 
 2 & 11298 &1.9 & 10399 & 1.9  \\
 4 & 5973 &3.7 & 5310 & 3.6  \\
 8 &  3502 & 6.4& 2881 &  6.7 \\
 16 & 2023  &11.0 & 1767  &  10.9 \\
 32 &   1685 & 13.2& 1424 &  13.5 \\
 64 &  1041 & 21.4& 880 &  21.9 \\
 128 &  571 &39.0 & 474 &  40.6 \\
 256 &  344 &64.7 & 287 &  67.1 \\
 512 &  213 &104.6 & 182 & 105.9 \\
 1024 & 149 &149.5 & 126 & 152.9
  \end{tabular}
\end{table}

\subsubsection{Weak scaling}

We next look into the weak scaling properties of the Poiseuille flow example and the GKB method with inner solver CG and fgmres (see Figure~\ref{fig:weak_prob3}). We do not present results for the two variants involving MUMPS, as they did not show any satisfactory weak scaling properties in our experiments. We scale the total workload with the increase in the core count such that each core is assigned 2048 ($=2\cdot 32^2$) elements. To always obtain an integer number of cores, we increase the total number of elements by 64 and 32 in $x$- and $y$-direction, respectively, up to the largest total problem of $n_x=2048, n_y=1024$ elements as in {\tt Prob 3}. We use the stopping tolerances $\tau=10^{-6}$ and $\tau_{in}=\{10^{-7}, 10^{-8}\}$, as required by {\tt Prob 3} for a precise computation, for all problem sizes. The number of iterations required for the GKB convergence is not constant across the obtained problem sizes. However, contrary to what is usually observed, it decreases monotonically with the increase in the problem size. Hence, the fluctuations noticeable in Fig.~\ref{fig:weak_prob3} may be due to algorithm-architecture interplay rather than an increased problem complexity when problem size grows. In particular, for smaller problem sizes, up to 400 cores (corresponding to meshes of up to $n_x=1280, n_y=640$ elements) these fluctuations are more pronounced and may be explained by the dominance of parallel overhead.
For the problem sizes beyond 400 cores the total time time stays about constant ( in the $\pm$10\% range).

\subsection{Isoviscous Stokes - 2D}\label{sec:q2p1}
Second, we adapt example {\tt ex62.c}  \footnote{\url{https://www.mcs.anl.gov/petsc/petsc-dev/src/snes/examples/tutorials/ex62.c.html}} given in the {\tt SNES} tutorials section in PETSc to comply with (\ref{eqn:stokes}). It simulates the 2D isoviscous variant of the Stokes problem (\ref{eqn:stokes}), approximated by a Q2-P1 
finite
element method with a discontinuous pressure field. 
The domain $\Omega$ is the unit square and is discretized by an unstructured mesh with quadrilateral elements. The exact solution is chosen as
\begin{align*}
\uv&=(x^3 +y^3, 2 x^3 - 3x^2y)\\
\pv&= \frac{3}{2} x^2 + \frac{3}{2} y^2-1. 
\end{align*}
Since we deal with a linear problem, the {\tt SNES} non-linear solver converges in one iteration. 
We focus thus on the solution of the system with the Jacobian matrix for which we use the {\tt PCFIELDSPLIT} environment with the GKB solver and compare the four solution techniques as in the previous example. 
We choose $n_x, n_y = 1024$ cells in the $x$- and $y$-direction, leading to $m\approx8.3 \cdot 10^6$ and $n\approx3.1 \cdot 10^6$ degrees of freedom. 
The finite element discretization errors (obtained by solution with MUMPS) are
\begin{align*}
\|\uv_h - \uv\|_{L^2} \approx 6.62\cdot 10^{-11}, \hspace{0.5cm} \|\pv_h - \pv\|_{L^2} \approx 1.51\cdot 10^{-7}. 
\end{align*}
As an implementation of the energy norm is not available in this example, we experimentally choose the stopping tolerances such that the solution respects the $L^2$-errors.  
We obtain $\tau=10^{-8}$ for the outer GKB iteration, and $\tau_{in}=10^{-13}$ for the inner iterative solvers in GKB-CG/fgmres. 
For GKB-MUMPS, we test several choices of $\nu$ (Table~\ref{tab:eta_stokes2d}) and choose $\nu=10^7$ for the following experiments. 
This large value is justified by the different scaling of the matrix blocks (contrary to the Poiseuille flow example, we did not apply central scaling to (\ref{eqn:saddlepoint})). 
The matrix $\Am \Am^T$ scales with $h^2\approx 10^{-6}$, and $\nu$ thus equilibrates the values in the blocks of (\ref{eqn:augsys_auglag}). 
Note that $\nu=10^8$ leads to an even smaller number of GKB iterations and a faster computation time. 
The ill-conditioning of the augmented matrix $\Mm = \Wm + \nu \Am \Am^T$ however leads to inaccuracies in the results, as can be noticed in a growing error $\mbox{err}_2^{\uv}$. 
Lastly, we choose $d=5$ in (\ref{eqn:stop}). 

\begin{table}
\caption{Choice of $\Nm = \frac{1}{\nu}$ for Stokes Q2-P1 in 2D. $\tau=10^{-8}$, $d=5$and 32 cores}
\begin{center}
 \begin{tabular}{r|r|r|r|r|r}
 \multicolumn{1}{c|}{$\nu$} & \multicolumn{1}{|c|}{$\mbox{err}_2^{\uv}$} & \multicolumn{1}{|c|}{$\mbox{err}_2^{\pv}$} & \multicolumn{1}{|c|}{l.b. estimate}  & \multicolumn{1}{|c|}{GKB iter} & \multicolumn{1}{|c}{time, $s$}\\
 \hline
0   & 7.23-11 & 1.51e-07  &   5.77e-09  &81 & 158 \\
$10^2$  & 7.52e-11 &  1.51e-07 &  9.91e-09 & 80 & 159\\
$10^5$  & 6.86e-11 &  1.51e-07 &   7.36e-09  & 66 & 150 \\
$10^7$  & 6.73e-11 & 1.51e-07  &  3.92e-09   & 18   & 116\\
$10^8$  & 1.21e-10 & 1.51e-07  &   5.60e-09  & 11   & 112\\
\end{tabular}
\end{center}
\label{tab:eta_stokes2d}
\end{table}

\subsubsection{Strong scaling}

The results are presented in Figure~\ref{fig:scal_q2p1} and Table~\ref{tab:scaleQ2P1}. 
As in the previous example, MUMPS and GKB-MUMPS are clearly the better choice for a small number of cores. For one core (and roughly for two cores), MUMPS and GKB-MUMPS are more than one order of magnitude faster than the iterative choices GKB-CG and GKB-fgmres. Compared to GKB-MUMPS, MUMPS on (\ref{eqn:augsys_auglag}) is faster for a small number of cores, but the two methods show a similar performance from about 32 cores onwards and reach a plateau at 64 cores. Once the computations pass from one to two nodes, i.e. from 32 to 64 cores, GKB-CG and GKB-fgmres start to outperform the methods involving MUMPS. For CG- and fgmres/BoomerAMG as inner solvers, the performance plateau is still not reached for 1024 cores. We observe a speed-up of about 517 for both GKB-CG and GKB-fgmres. The fastest computation time of 8 seconds is reached for GKB with the inner solver fgmres preconditioned with BoomerAMG.

\begin{table}
 \caption{Solver time and strong scaling for Stokes 2D with Q2-P1 example. } 
 \begin{center}
 \begin{tabular}{l||l|l||l|l|p{1.6cm}|p{1.6cm}|p{1.6cm}|p{1.5cm}}
  &\multicolumn{2}{c||}{MUMPS}  & \multicolumn{2}{|c|}{GKB-MUMPS} & \multicolumn{2}{c|}{GKB-CG} &\multicolumn{2}{c}{GKB-fgmres}\\
 & \multicolumn{2}{c||}{}&  \multicolumn{2}{c|}{\scriptsize{($\nu$=$10^7$), $\tau$=$10^{-6}$}} & \multicolumn{2}{c|}{\scriptsize{($\nu$=0,$\tau$=$10^{-8}$,$\tau_{in}$=$10^{-13}$)}} &\multicolumn{2}{c}{\scriptsize{($\nu$=$0$,$\tau$=$10^{-8}$,$\tau_{in}$=$10^{-13}$)}}\\
 cores & time (s) & scale & time (s) & scale & time (s) & scale& time (s)& scale\\ 
  \hline
 1 & 223 & &  328 & &  5174  & & 4133 & \\
 2 & 182 & 1.2&  232 & 1.4 & 2337 &2.2 &  1963 & 2.1  \\
 4 & 144 & 1.5 &  165 &2.0 & 958 &5.4 & 778 & 5.3  \\
 8 & 125 & 1.8 &  136 &2.4 & 563 &9.2 & 434 &  9.5 \\
 16 & 115 & 1.9 & 119 &2.8 & 314 &16.5 & 246 &  16.8 \\
 32 &112  & 2.0 & 114 &2.9 & 226 &22.9 & 168 &  24.6 \\
 64 & 108  & 2.1 & 107 & 3.1 & 115 & 45.0& 85 &  48.6 \\
 128 & 107 & 2.1 & 107 &3.1 & 56 &92.4 & 44 &  93.9 \\
 256 &  109& 2.0 & 107 &3.1 & 31 &166.9 & 23 &  179.7 \\
 512 & 110  & 2.0 & 109 & 3.0& 17 &304.4 & 13 & 317.9 \\
 1024 &  118 & 1.9 & 129 & 2.6& 10 &517.4 & 8 & 516.6
  \end{tabular}
 \end{center} 
  \label{tab:scaleQ2P1}
 \end{table}

\subsubsection{Weak scaling}
\begin{figure}
\centering
\begin{subfigure}{0.45\textwidth}
\includegraphics[width=\textwidth]{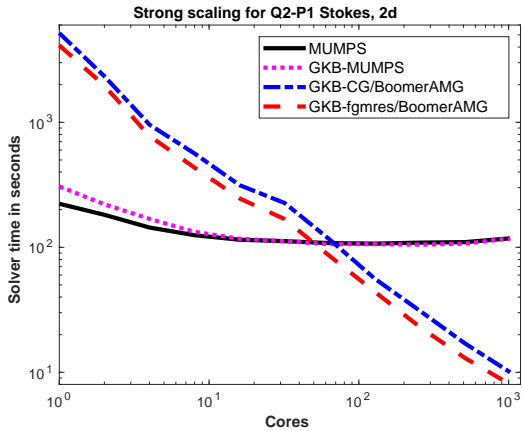}
\caption{Strong scaling 2D isoviscous Stokes.}
\label{fig:scal_q2p1}
\end{subfigure}
\begin{subfigure}{0.45\textwidth}
\includegraphics[width=\textwidth]{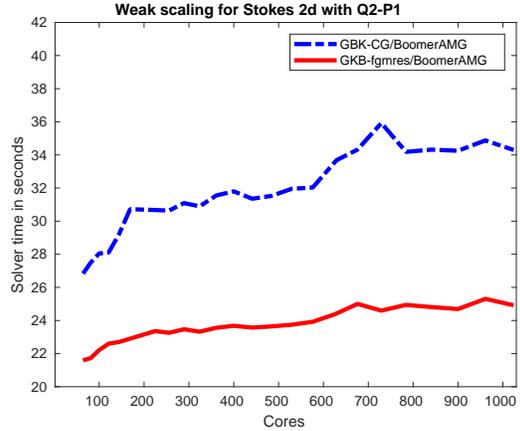}
\caption{Weak scaling 2D isoviscous Stokes.}
\label{fig:weak_st_2d}
\end{subfigure}
\caption{Scaling behavior on Kraken.}
\label{fig:stokes_2d}
\end{figure}

We next look into the weak scaling properties of the GKB algorithm for the 2D Stokes Q2-P1 example (Fig.~\ref{fig:weak_st_2d}). We distribute the workload such that each core obtains 4096 ($=64^2$) elements. For this, we increase the mesh size in steps of 64 in $x$- and $y$-direction up to $n_x=2048, n_y=2048$ and compute the required (integer) number of cores. We choose $\tau=10^{-8}$ and $\tau_{in}=10^{-13}$ as stopping tolerances, such that the largest problem is solved with sufficient precision (the required tolerances were found experimentally). Although the total problem size increases, the number of GKB iterations stays between 80 and 82 for any mesh size. The total number of iterations until convergence, i.e. outer $\times$ inner iterations, increases slightly for the CG method towards the end. In general, the rather flat timing curves for both methods, presented in Fig.~\ref{fig:weak_st_2d}, point to a good weak scaling behavior.     

\subsection{Comparison with AMD architecture}

To test the behavior of the GKB algorithm on another architecture, we ran the 2D Stokes Q2-P1 test on two bi-socket AMD Rome (EPYC 7702) nodes at 2 Ghz. Each socket contains 64 cores and shares 256 GB DDR4 memory. We compare two configurations. The first configuration (config 1) uses 4 NUMA (non-uniform memory access) nodes per socket (NPS4). The second configuration (config 2) uses 1 NUMA node per socket (NPS1), which corresponds to the memory access configuration of the Intel Skylake sockets on Kraken. The results are given in Table~\ref{tab:scaleStokesAMD}. Up to a parallelism of 4 cores, the computations on both AMD configurations are faster in absolute computation time than the computations on the Intel nodes, although having a smaller clock rate (2 Ghz compared to 2.5 Ghz).  For config 1, we also observe a higher scalability up to 64 nodes, when it slows down and the absolute computation time becomes higher. The final speed-up at 128 cores is about 40 (compared to 92 on Intel) for GKB-CG and 43 for GKB-fgmres (compared to 94 on Intel). Although having the fastest sequential performance, config 2 is the least interesting set-up in this comparison, as the scalability already slows down significantly at 6 cores. It obtains a final speed-up of about 30 for GKB-CG and 33 for GKB-fgmres. 

\begin{table}
 \caption{Solver time and strong scaling for 2D Stokes Q2-P1 on AMD architecture } 
 \begin{center}
\begin{tabular}{l||p{1.3cm}|p{1.3cm}|p{1.5cm}|p{1.3cm}|p{1.3cm}|p{1.3cm}|p{1.3cm}|p{1.3cm}}
  & \multicolumn{2}{c|}{CG - config 1} &\multicolumn{2}{c|}{CG - config 2}& \multicolumn{2}{c|}{fgmres - config 1} &\multicolumn{2}{c}{fgmres - config 2}\\
 & \multicolumn{2}{c|}{\scriptsize{($\nu$=0,$\tau$=$10^{-8}$,$\tau_{in}$=$10^{-13}$)}} &\multicolumn{2}{c|}{\scriptsize{($\nu$=$0$,$\tau$=$10^{-8}$,$\tau_{in}$=$10^{-13}$)}} & \multicolumn{2}{c|}{\scriptsize{($\nu$=0,$\tau$=$10^{-8}$,$\tau_{in}$=$10^{-13}$)}} &\multicolumn{2}{c}{\scriptsize{($\nu$=$0$,$\tau$=$10^{-8}$,$\tau_{in}$=$10^{-13}$)}}\\
 cores & time (s) & scale& time (s)& scale& time (s) & scale& time (s)& scale\\ 
  \hline
   1 &4268 & & 3672 &   & 3552   & &3072&\\
 2 & 1721 & 2.5& 1741 & 2.1  &1439 & 2.5&1450 &2.1\\
 4 & 767 &5.6 & 852 & 4.3 & 615&5.8 &688 &4.5\\
 8 & 429  &9.9 & 632 & 5.8 &  331 &10.7 & 482 &6.4\\
 16 & 245  &17.4 & 523  & 7.0&188 &18.9 & 402 &7.6\\
 32 & 173  &24.7 &284  & 12.9  &127 &28.0 &210 &14.6\\
 64 &  145 &29.4 & 162 &  22.7 &107 &33.2 & 119&25.8\\
 128 & 108  & 39.5& 122 & 30.1  & 83 & 42.8 & 94 &32.7\\
   \end{tabular}
  \end{center}
  \label{tab:scaleStokesAMD}
 \end{table}

\subsection{Isoviscous Stokes - 3D}
As last example, we consider a 3D isoviscous Stokes problem and adapt {\tt ex62.c} in the {\tt SNES} tutorials section in PETSc to comply for use with our solver. Again, we use a Q2-P1 finite element discretization with a discontinuous pressure field. The domain $\Omega$ is the unit cube and it is discretized by an unstructured mesh. The exact solutions used are 
\begin{align*}
\uv &= (x^3 + y^3 + 3 z^2 x, 2 x^3 - 3 x^2 y, -z^3 x + y^2)\\
\pv &= 1.5 x^2 + 1.5 y^2 + 1.5 z^2 - 1.5. 
\end{align*}
For the following strong scaling test, we use a mesh with $n_x, n_y, n_z = 64$ elements, which leads to $m \approx 7.2\cdot 10^6$ and degrees of freedom for the (1,1)-block and $n\approx 6.1\cdot 10^6$ constraints. Contrary to the 2D problem, the ratio between the number of physical degrees of freedom and constraints increased and is roughly $\frac{7}{6}$. The finite element discretization errors, obtained approximately by using GKB-MUMPS with a small tolerance $\tau$, are $\|\uv_h - \uv\|_{L^2} \approx 3.11\cdot 10^{-7}$ and $\|\pv_h - \pv\|_{L^2} \approx 4.72\cdot 10^{-5}$. Although the problem is linear, the non-linear solver needs a rather low tolerance for the GKB method, i.e. to solve the system defined by the Jacobi matrix, to converge in one step. We thus choose $\tau = 10^{-10}$ for the outer GKB iteration and $\tau_{in}=10^{-11}$ for the inner iterative solvers CG or fgmres. The parameter in the augmented Lagrangian approach is chosen as $\nu=0$ and for the stopping criterion we use $d=5$. The matrices for 3D finite element problems are usually much denser than their 2D counterparts. As a consequence, the numerical simulations require a significant amount of memory. In our experiments, we use 20 of 38 available cores on each node, as the algorithm suffered of memory contention for a higher number of cores per nodes.

\subsection{Strong scaling}
\begin{figure}
\centering
\begin{subfigure}{0.45\textwidth}
\includegraphics[width=\textwidth]{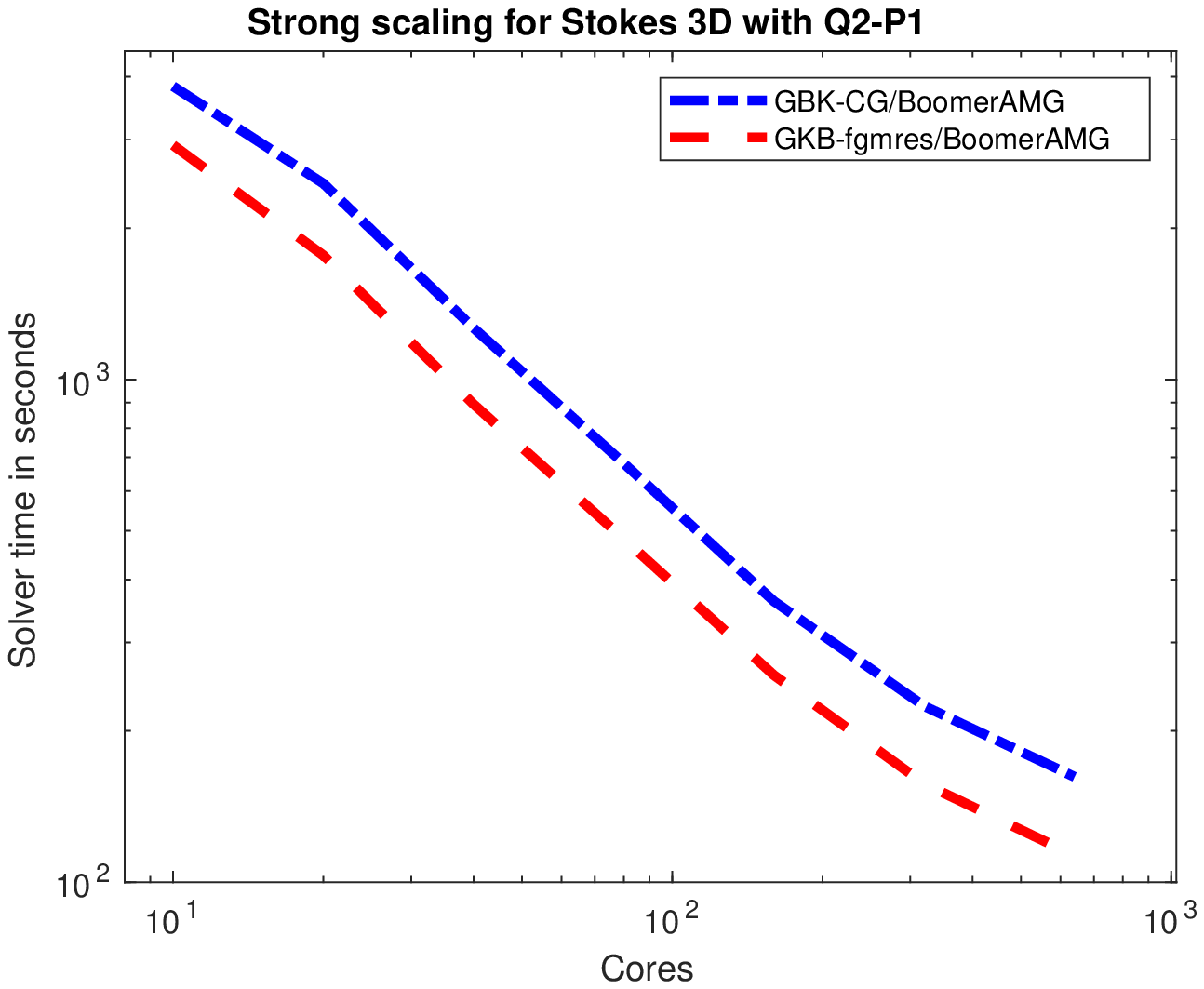}
\caption{Strong scaling 3D isoviscous Stokes.}
\label{fig:str_ex62_3d}
\end{subfigure}
\begin{subfigure}{0.45\textwidth}
\includegraphics[width=\textwidth]{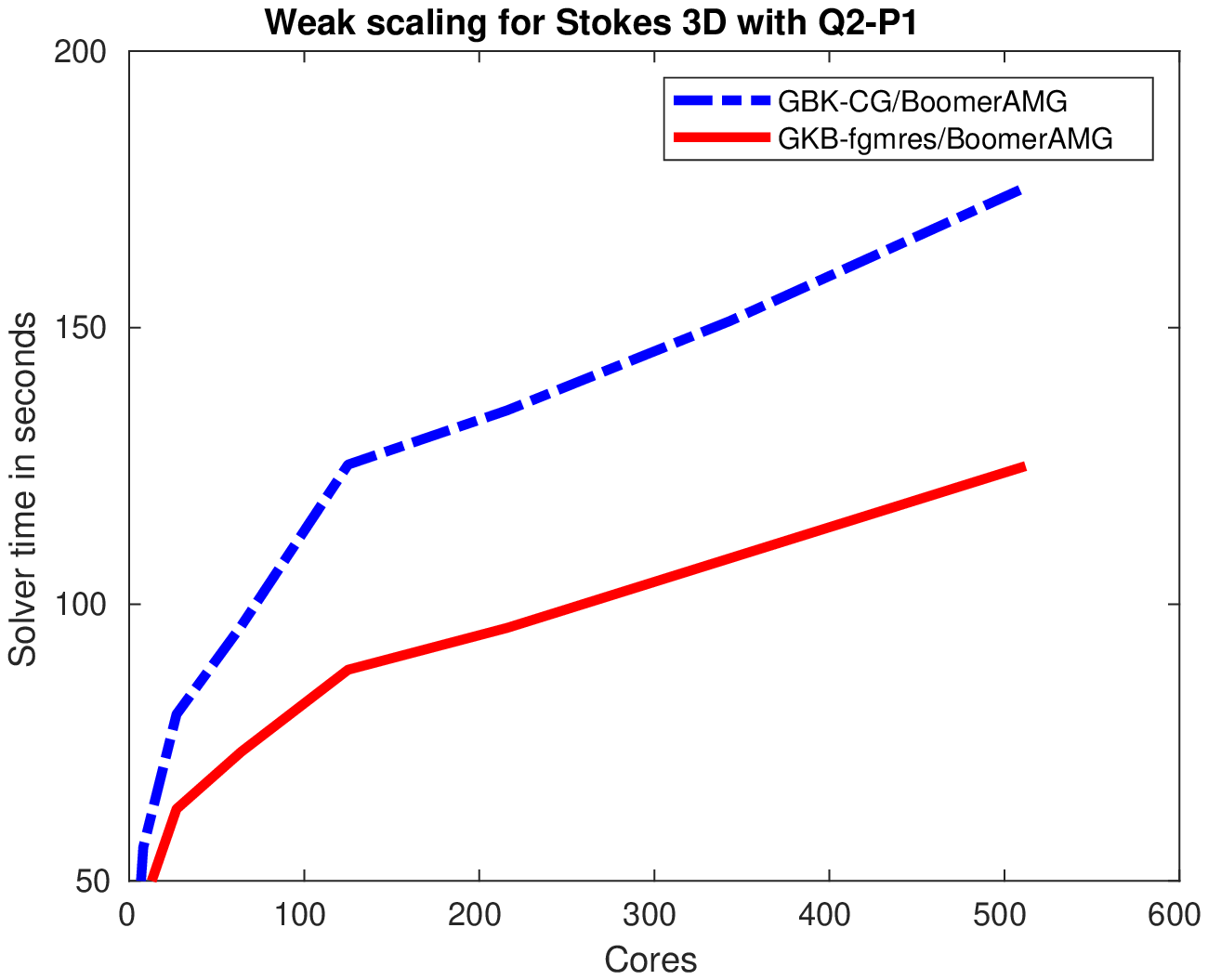}
\caption{Weak scaling 3D isoviscous Stokes.}
\label{fig:weak_ex62_3d}
\end{subfigure}
\caption{Scaling behavior on Kraken.}
\end{figure}

When using only a low number of cores per node on Kraken, it was in principle possible to run the GKB-MUMPS algorithm. The computation times were, however, not competitive. We thus compare the strong scaling results for GKB-CG and GKB-fgmres only. The results are given in Table~\ref{tab:scaleStr3d} and Fig.~\ref{fig:str_ex62_3d}. As in the 2D case, the GKB method with inner iterative solvers scales for a growing number of cores. For both methods, the plateau is still not reached for 640 cores. GKB-CG obtains a speed-up of about 24 (with optimal value of 64 in our scale) and GKB-fgmres reaches a speed-up of about 26. In terms of the total computation time, GKB-fgmres is the faster choice for any number of cores. As we have also observed in the two-dimensional examples before, the CG solver needs on average one iteration more than fgmres does so for convergence in each GKB iteration. 

  \begin{table}
  \caption{Solver time and strong scaling for Stokes 3D.}
    \centering
    \begin{tabular}{l||p{1.4cm}|p{1.4cm}|p{1.4cm}|p{1.4cm}}
  & \multicolumn{2}{c|}{GKB-CG} &\multicolumn{2}{c}{GKB-fgmres}\\
 & \multicolumn{2}{c|}{\scriptsize{($\nu$=0,$\tau$=$10^{-10}$,$\tau_{in}$=$10^{-11}$)}} &\multicolumn{2}{c}{\scriptsize{($\nu$=$0$,$\tau$=$10^{-10}$,$\tau_{in}$=$10^{-11}$)}}\\
 cores & time (s) & scale& time (s)& scale\\ 
  \hline
 10 & 3831  &  & 2925  &  \\ 
 20 & 2454 &1.6 & 1766 & 1.7  \\
 40 & 1266 &3.0 & 893 &  3.3 \\
 80 & 681  &5.6 & 482 &  6.1 \\
 160 & 362  &10.6 & 258  & 11.3  \\
 320 & 224  &17.1 & 157 &  18.3 \\
 640 & 162 & 23.6& 113 &  25.9 \\
  \end{tabular}
 \label{tab:scaleStr3d}
\end{table} 
 
    
\subsection{Weak scaling}

For the weak scaling results of the 3D Stokes Q2-P1 example, we keep constant the 512 (=$8^3$) elements per core. We, therefore, increase the total problem size by 8 in $x$-, $y$-,  and $z$-directions, from which we compute the required (integer) number of cores. The results are presented in Table~\ref{tab:WeakSc3d} and Fig.~\ref{fig:weak_ex62_3d}. The total number of iterations, i.e. outer GKB iterations $\cdot$ inner iterations, increases only little for CG and stays constant for the fgmres inner solver from 64 cores onwards. In terms of computation times, we observe for larger core numbers (125+ cores) in Fig.~\ref{fig:weak_ex62_3d}, similar to  Fig.~\ref{fig:weak_st_2d},  a linear growth. The steeper slope of the timing curves may be explained by the memory wall starting to affect GKB with iterative inner solvers for large 3D problems and not enough memory per node with the maximum of 38 cores used in our experiments. 

 \begin{table}
 \caption{Solver time and weak scaling for Stokes 3D.}
 \centering
    \begin{tabular}{l|l||p{1.4cm}|p{1.4cm}|p{1.4cm}|p{1.4cm}}
  & & \multicolumn{2}{c|}{GKB-CG} &\multicolumn{2}{c}{GKB-fgmres}\\
 && \multicolumn{2}{c|}{\scriptsize{($\nu$=0,$\tau$=$10^{-10}$,$\tau_{in}$=$10^{-11}$)}} &\multicolumn{2}{c}{\scriptsize{($\nu$=$0$,$\tau$=$10^{-10}$,$\tau_{in}$=$10^{-11}$)}}\\
 cores & $n$ & time (s) & total iter & time (s)& total iter \\ 
  \hline
 27 & 24 & 80 &1809& 63 &   1556\\
 64 & 32 & 96& 1894& 73  & 1622 \\
 125 & 40 &125 &2138 & 88 & 1622 \\
 216 & 48 &135 &2144 & 96 & 1622  \\
 343 & 56 &151 & 2158& 108 & 1622  \\
 512 & 64 & 175 & 2218&  125  & 1622\\
  \end{tabular}
\label{tab:WeakSc3d}
\end{table}

\section{Conclusions}
We presented an iterative algorithm based on the Golub Kahan bidiagonalization for 2x2 block matrices. Furthermore, we outlined our PETSc implementation of the solver and applied it to the Poiseuille flow as well as an 2D and 3D isoviscous Stokes problem. The strong scaling results showed that MUMPS and the GKB-MUMPS solvers start to level off at about 64 cores for a constant problem size, with speed-up factors of at most 2.9, while the GKB algorithm with BoomerAMG as preconditioner showed a speed-up between 40 to 517 at 1024 cores. Regarding the gains in  the absolute computation time, either MUMPS or GKB-MUMPS are the methods of choice for up to 32 cores and unknowns of the order of $10^6$. When more than 64 cores are available, the GKB method with inner iterative solvers outperforms its direct counterparts, and hence, should be employed. When increasing the problem size to about $2\cdot 10^{7}$ unknowns ({\tt Prob 4}), MUMPS required more memory than that available on Cerfacs' cluster Kraken in either the direct or GKB-MUMPS case. Iterative methods, such  as GKB-fgmres and GKB-CG, are usable for this problem and present the only alternative to its solution.
In addition, they scale well with the increase in total problem size and the cores counts when the workload (the number of elements) is kept constant per core. 

\section*{Acknowledgements} The work of the second author was supported in part by the U.S. Department of Defense High Performance
Computing Modernization Program, through a HASI grant, and by the National
Science Foundation under grant CNS-1828593.

\end{document}